\documentclass[11pt]{article}
\usepackage[T1]{fontenc}
\usepackage[a4paper]{anysize}\marginsize{2.0cm}{2.0cm}{1.0cm}{1.0cm}
\pdfpagewidth=\paperwidth \pdfpageheight=\paperheight
\usepackage{pgf,tikz,hyperref}
\usetikzlibrary{arrows}
\usepackage{amsfonts,amssymb,amsthm,amsmath,eucal}
\usepackage{graphicx}
\usepackage{ltablex}
\usepackage{caption, booktabs}
\usepackage{longtable}
\usepackage{pgfplots}
\pgfplotsset{compat=1.15}
\usepackage{mathrsfs}

\theoremstyle{plain}
\newtheorem{thm}{Theorem}[section]
\newtheorem{theorem}[thm]{Theorem}

\newtheorem{lemma}[thm]{Lemma}
\newtheorem{proposition}[thm]{Proposition}

\theoremstyle{definition}
\newtheorem{definition}[thm]{Definition}
\newtheorem{remark}[thm]{Remark}

\newtheorem{notation}[thm]{Notation}

\newtheorem{thevarthm}[thm]{\varthmname}

\newenvironment{varthm*}[1]{\trivlist\item[]{\bf #1.}\it}{\endtrivlist}

\newcommand\be{\begin{eqnarray*}}
\newcommand\ee{\end{eqnarray*}}

\newcommand\cald{\mathcal D}

\renewcommand\P{\mathbb P}

\newcommand\newop[2]{\def#1{\mathop{\rm #2}\nolimits}}
\newop\edim{edim}
\newop\Zeroes{Zeroes}
\newop\Jac{Jac}
\newop\Ass{Ass}
\newop\SL{SL}
\newop\PGL{{\P}GL}
\newop\Km{Km}
\newop\reg{reg}

\newcolumntype{L}{>{$}l<{$}}

\newcommand\keywords[1]{{\renewcommand\thefootnote{}\footnotetext{\textit{Keywords:} #1.}}}
\newcommand\subclass[1]{{\renewcommand\thefootnote{}\footnotetext{\textit{Mathematics Subject Classification (2010):} #1.}}}

\definecolor{zzttff}{rgb}{0.6,0.2,1}
\definecolor{zzttqq}{rgb}{0.6,0.2,0}
\definecolor{qqqqff}{rgb}{0,0,1}
\definecolor{qqwuqq}{rgb}{0,0.39215686274509803,0}

\begin{document}

\author{\L{}ukasz Merta \and Marcin Zieli\'nski}
\title{On quartics with the maximal number of the maximal tangency lines}
\date{\today}
\maketitle
\thispagestyle{empty}

\begin{abstract}
In this note, we examine the arrangements of lines and configurations of points that emerge from Fermat (von Dyck) and Komiya-Kuribayashi quartics. These quartics are characterized by having the maximum number of lines of maximal tangency, that is, lines for which the intersection multiplicity at the tangency point is equal to the degree of the curve. Additionally, we delve into the study of sextactic points on these quartics -- points at which there exists a conic with the curve having a local intersection multiplicity of at least 6, which is one more than that observed at a general point -- alongside the related configurations of conics.
\end{abstract}
\keywords{quartic curves, bitangents}
\subclass{MSC 14C20 \and MSC 14N20}

\section{Introduction}
It is well known that there are exactly $28$ bitangents to a smooth plane quartic curve (in fact, around 1839 Pl\"ucker established that a smooth plane curve of degree $d$ has $\frac{1}{2}d(d-2)(d^2-9)$ bitangents, see \cite{GH}). Sometimes it happens that the two points where the tangent line touches the quartic fall together to a point, where the order of tangency is maximal, hence $4$. We call these lines the maximal tangency lines. Such lines can be defined for an arbitrary plane curve.
\begin{definition}
    Let $C$ be a smooth plane curve of degree $d$. We say that $L$ is a \emph{maximal tangency line} (MTL for short) of $C$, if it intersects $C$ in exactly one point $P$. We say then that $P$ is the \emph{maximal tangency point} (MTP for short) of $C$.
\end{definition}
If $C$ is a line, then any line in the plane different than $C$ is its MTL. If $C$ is a conic, then any line tangent to $C$ is its MTL.
In the case of cubics the maximal tangency lines are the tangents in the flex points. Their configuration and the configuration dual to the flex points have interesting properties from the combinatorial and algebraic point of view, see e.g. \cite{DST2013}. This motivates our interest in the arrangements related to the maximal tangency lines for quartics. Contrary to the case of cubics, where there are always exactly $9$ maximal tangency lines, for quartics their number ranges from $0$ to $12$. The points of tangency are called in this case \emph{undulation} points, see \cite{PopSha15}.

Here we focus on quartics with the maximal number of MTLs. It was shown in \cite{KK1979} that there are exactly two quartics with $12$ maximal tangency lines, namely the Fermat (or von Dyck) quartic $F$ with the equation
$$F:\; x^4+y^4+z^4=0$$
and the Komiya-Kuribayashi quartic $K$ given by the equation
$$K:\; x^4 + y^4 + z^4 + 3(x^2y^2 + x^2z^2 + y^2z^2)=0.$$
Both quartics have relatively large automorphism groups: of order $96$ for the Fermat quartic and $24$ for the Komiya-Kuribayashi quartic. Note that the order of the automorphism group of the Klein quartic is even bigger, namely $168$, which is the maximal possible order of the automorphism group of a curve of genus $3$ by the well known Hurwitz bound, see \cite{Forster81}. However the Klein quartic has no MTLs, see \cite{KK1979}.

Recently Roulleau \cite{Xavier2} proposed a systematic way to construct new line arrangements and point configurations out of given ones. In particular he introduced a family of operators providing a uniform way to run the constructions. We recall some of these operators below.


\begin{definition}{(Roulleau)}\label{roulleau_L}
    For a subset ${\bf n} \subset \mathbb{N} \setminus \{0,1\}$ and an arrangement of lines $L$ in the plane, we define
    $\mathcal{P}_{\bf n}(L)$ as a set of points where exactly $k$ lines from $L$ intersect, for some $k \in {\bf n}$.
\end{definition}

To simplify the notation, in the case ${\bf n} = \{k\}$ for some $k \in \mathbb{N}$, we will write $\mathcal{P}_k$ instead of $\mathcal{P}_{\{k\}}$.

\begin{definition}{(Roulleau)}\label{roulleau_D}
    We define $\mathcal{D}$ as an operator which maps a point $(a : b : c) \in \mathbb{P}^2$ to a line given by $ax + by + cz = 0$. We also define an operator $\check{\mathcal{D}}$ as the inverse of $\mathcal{D}$.
\end{definition}

These dual operators mentioned above can also be applied to configurations of points and lines, respectively.

Here we also recall merely the following fundamental combinatorial invariant of line arrangements.
\begin{notation}
    Let ${\bf L}=\left\{L_1,\ldots,L_d\right\}$ be an arrangement of $d$ lines. Let $r$ be the maximal multiplicity of ${\bf L}$ which is defined as 
    $r = {\rm max} \{k : \text{ arrangement } {\bf L} \text{ admits } k\text{-fold points}\}$.
    The $t$-vector of ${\bf L}$ is the vector
    $$t({\bf L})=(t_2,t_3,\ldots,t_r),$$
    where $t_i$ is the number of points where exactly $i$ lines from ${\bf L}$ intersect.
\end{notation}
Let us now briefly outline the content of this paper. In Section \ref{sec2} we study particular line arrangements which are related to the MTPs and the MTLs to the Fermat quartic. In Section \ref{sec3} we describe some more line arrangements and their properties, but this time it is connected to the MTPs and the MTLs to the Komiya-Kuribayashi quartic. Finally, in section \ref{sec4} we compute the exact coordinates of sextactic points on both quartics, using a method involving the second Hessian, while in Section \ref{sec5}, we study the arrangements of conics which are tangent to a given quartic at two sextactic points with multiplicity $4$ at each such point. Our study in this direction is inspired by a recent work of Szemberg and Szpond \cite{SzeSzp24}, which was later continued in a paper by Merta and Zięba \cite{MERTA}. 


\section{The Fermat quartic}\label{sec2}
For the Fermat quartic, it is easy to check that the $12$ maximal tangency lines $LF = \{ LF_1,\ldots, LF_{12}\}$ are linear factors of the polynomial
$$H=(x^4+y^4)(y^4+z^4)(z^4+x^4).$$

\begin{lemma}
The $12$ maximal tangency lines $LF_1, \dots, LF_{12}$ intersect in $48$ double points and $3$ points of multiplicity $4$.
\end{lemma}

\begin{proof}
    Necessary easy computations can be done in principle by hand.
    We utilized the Singular symbolic algebra program \cite{Singular}.
\end{proof}

In particular we obtain in that case the $t$-vector:
$$t=(48,0,3).$$
The points of multiplicity $4$ are of course common intersection points of lines determined by linear factors of $x^4+y^4$ (respectively: $y^4+z^4$ and $z^4+x^4$), so that these are the coordinate points $(0:0:1), (1:0:0)$ and $(0:1:0)$. The remaining $48$ points, i.e.\ $\mathcal{P}_2(LF)$, are determined by pairs of lines passing through \emph{distinct} coordinate points.


\begin{proposition}
    Let $\mathcal{I}$ be the saturated ideal of the $48$ points determined by $\mathcal{P}_2(LF)$. Then $\mathcal{I}$ is generated in degree $8$ and its minimal free resolution has the following shape:
    $$0 \longrightarrow S(-12)^2 \longrightarrow S(-8)^3 \longrightarrow S \longrightarrow S/\mathcal{I} \longrightarrow 0.$$
\end{proposition}

The generators of $\mathcal{I}$ can be written in the following form:
\begin{equation}\label{ideal_4}(x^4 + y^4)(x^4 + z^4),\ (x^4 + y^4)(y^4 + z^4),\ (x^4 + z^4)(y^4 + z^4).
\end{equation}
Interestingly, the $12$ points of tangency $PF = \{PF_1,\ldots,PF_{12}\}$ are the duals of these lines, i.e.\ $PF = \check{\mathcal{D}}(LF)$. The ideal of those points is a complete intersection generated by
$$xyz \mbox{ and } x^4+y^4+z^4,$$
see our Singular code and also \cite{DYCK}. Thus the points are distributed by $4$ on the $3$ coordinate axes. 
\begin{theorem}
    The $4$ MTPs of the Fermat quartic sitting on each of the coordinate axes built a harmonic four.
\end{theorem}
\begin{proof}
    By symmetry it suffices to show it for one axis only. Let $u$ be a primitive root of unity of order $8$. Then the relevant points on the $z=0$-axis are
    $$(1:u:0),\; (1:u^3:0),\; (1:u^5:0)\mbox{ and } (1:u^7:0).$$
    The cross-ratio of these points will not change when we divide the second coordinate by $u$ (because it corresponds to the projective automorphism of $\P^2$ defined by 
    $$(x:y:z)\longleftrightarrow (x:u^7\cdot y:z).$$
    In affine coordinates on the $y$-axis the points correspond then to $1,i,-1$ and $-i$ and it is well known that they are a harmonic four.
\end{proof}

\begin{figure}[ht!]
\centering
\begin{tikzpicture}[line cap=round,line join=round,>=triangle 45,x=1cm,y=1cm]
\clip(-7.58,-5.04) rectangle (7.58,5.88);
\draw [line width=2pt,domain=-7.58:7.58] plot(\x,{(-30-0*\x)/10});
\draw [line width=2pt,domain=-7.58:7.58] plot(\x,{(-20--7*\x)/-5});
\draw [line width=2pt,domain=-7.58:7.58] plot(\x,{(-20-7*\x)/-5});
\draw [line width=0.8pt,color=qqwuqq,domain=-7.58:7.58] plot(\x,{(--11.116216216216216--5.968918918918918*\x)/2.263513513513513});
\draw [line width=0.8pt,color=qqwuqq,domain=-7.58:7.58] plot(\x,{(--10.697297297297297-4.568918918918919*\x)/4.736486486486487});
\draw [line width=0.8pt,color=qqwuqq,domain=-7.58:7.58] plot(\x,{(--16.8-1.4*\x)/-7});
\draw [line width=0.8pt,color=qqqqff,domain=-7.58:7.58] plot(\x,{(-8.918054054054055--2.902162162162162*\x)/4.327027027027027});
\draw [line width=0.8pt,color=qqqqff,domain=-7.58:7.58] plot(\x,{(-5.323877282688093--1.9675675675675675*\x)/-4.448648648648649});
\draw [line width=0.8pt,color=qqqqff,domain=-7.58:7.58] plot(\x,{(-7.1824864864864875-4.8697297297297295*\x)/0.12162162162162149});
\draw [line width=0.8pt,color=zzttqq,domain=-7.58:7.58] plot(\x,{(--3.5855441928414873--1.5645945945945945*\x)/4.92837837837838});
\draw [line width=0.8pt,color=zzttqq,domain=-7.58:7.58] plot(\x,{(--7.04864864864865-4.332432432432432*\x)/-0.9054054054054061});
\draw [line width=0.8pt,color=zzttqq,domain=-7.58:7.58] plot(\x,{(--9.301081081081083--2.7678378378378374*\x)/-4.022972972972973});
\draw [line width=0.8pt,color=zzttff,domain=-7.58:7.58] plot(\x,{(--10.55193571950329--4.654054054054054*\x)/4.891891891891891});
\draw [line width=0.8pt,color=zzttff,domain=-7.58:7.58] plot(\x,{(--10.246702702702702-5.902702702702703*\x)/1.9362162162162173});
\draw [line width=0.8pt,color=zzttff,domain=-7.58:7.58] plot(\x,{(--17.088--1.2486486486486488*\x)/-6.8281081081081085});
\draw (-7.04,-2.39) node[anchor=north west] {$z = 0$};
\draw (4.6,-4.02) node[anchor=north west] {$y = 0$};
\draw (0.94,5.48) node[anchor=north west] {$x = 0$};
\draw (-1.38,3.15) node[anchor=north west] {$u$};
\draw (-3.08,-3.06) node[anchor=north west] {$u$};
\draw (3.65,-1.72) node[anchor=north west] {$u$};
\draw (-2.35,2.04) node[anchor=north west] {$u^3$};
\draw (2.71,0.69) node[anchor=north west] {$u^3$};
\draw (-1.41,-2.99) node[anchor=north west] {$u^3$};
\draw (-3.34,0.6) node[anchor=north west] {$u^5$};
\draw (-4.5,-1) node[anchor=north west] {$u^7$};
\draw (0.31,-2.99) node[anchor=north west] {$u^5$};
\draw (1.98,2.12) node[anchor=north west] {$u^5$};
\draw (2.14,-2.99) node[anchor=north west] {$u^7$};
\draw (0.92,3.29) node[anchor=north west] {$u^7$};
\begin{scriptsize}
\draw [fill=black] (-4.108108108108108,-1.7513513513513512) circle (2.5pt);
\draw [fill=black] (-3.022972972972973,-0.23216216216216257) circle (2.5pt);
\draw [fill=black] (-1.5216216216216218,1.8697297297297295) circle (2.5pt);
\draw [fill=black] (-0.7364864864864866,2.9689189189189187) circle (2.5pt);
\draw [fill=black] (0.7837837837837833,2.902702702702703) circle (2.5pt);
\draw [fill=black] (1.9054054054054061,1.332432432432432) circle (2.5pt);
\draw [fill=black] (2.927027027027027,-0.09783783783783795) circle (2.5pt);
\draw [fill=black] (4,-1.6) circle (2.5pt);
\draw [fill=black] (2.72,-3) circle (2.5pt);
\draw [fill=black] (1,-3) circle (2.5pt);
\draw [fill=black] (-1.4,-3) circle (2.5pt);
\draw [fill=black] (-3,-3) circle (2.5pt);
\end{scriptsize}
\end{tikzpicture}
\caption{Configuration of $12$ lines passing through the MTPs on the Fermat quartic $x^4 + y^4 + z^4$.}
\label{fig:arrangement}
\end{figure}
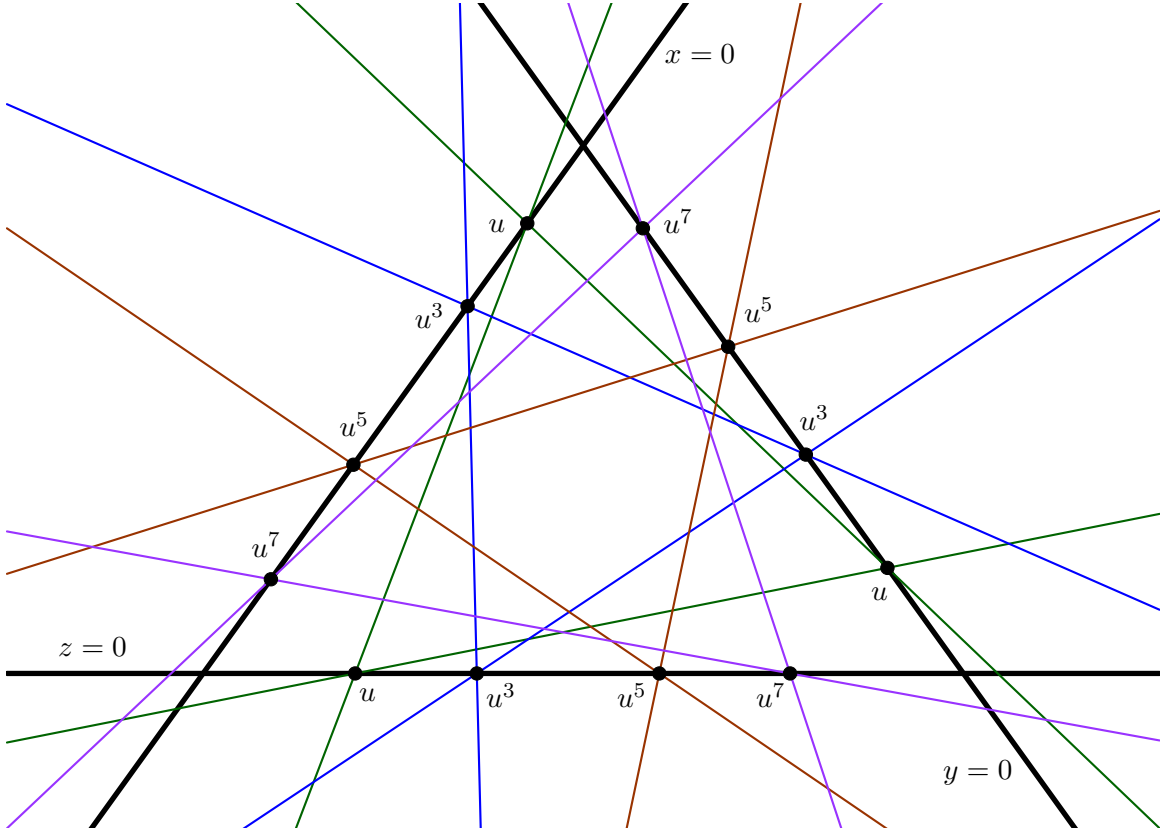

The coordinates of all $12$ MTPs can be written as $(0 : 1 : u^k),\ (u^k : 0 : 1)$ and $(1 : u^k : 0)$ for $k \in \{1, 3, 5, 7\}.$ For each such $k$, these three points lie on three lines (2 points per line), namely
\begin{align*}
    x + u^{2k}y - u^kz &= 0, \\
    y + u^{2k}z - u^kx &= 0, \\
    z + u^{2k}x - u^ky &= 0.
\end{align*}
Gathering these lines for all $k$, we obtain a symmetric arrangement of $12$ lines. The simplified illustration of this arrangement is presented in Figure \ref{fig:arrangement}. However, it turns out that these lines intersect in $66$ distinct points ($54$ points besides the MTPs) and there is no intersection point with multiplicity higher than $2$.

\section{The Komiya-Kuribayashi quartic}\label{sec3}
In this case the equations of maximal tangency lines were computed explicitly by Edge \cite{Edge1945}:\\
$$\begin{array}{llll}
LK_1:\;-2ix-y+z,& LK_2:\;-2ix+y+z,& LK_3:\;2ix-y+z,& LK_4:\;2ix+y+z,\\
LK_5:\;-2iy-x+z,& LK_6:\;-2iy+x+z,& LK_7:\;2iy-x+z,& LK_8:\;2iy+x+z,\\
LK_9:\;-2iz-y+x,& LK_{10}:\;-2iz+y+x,& LK_{11}:\;2iz-y+x,& LK_{12}:\;2iz+y+x.
\end{array}$$

These lines are general with respect to their intersection, i.e., they intersect by two in altogether 66 distinct points. These points in turn are far from being general. Direct inspection by Singular \cite{Singular} shows that their ideal has linear resolution.

\begin{proposition}
    Let $\mathcal{J}$ be the saturated ideal of the $66$ points determined by $\mathcal{P}_2(LK)$. Then $\mathcal{J}$ is generated in degree $11$ and its minimal free resolution has the following shape:
    $$0 \longrightarrow S(-12)^{11} \xrightarrow{\mbox{ }\alpha\mbox{ }} S(-11)^{12} \longrightarrow S \longrightarrow S/\mathcal{J} \longrightarrow 0.$$
\end{proposition}

\begin{proof}
    The resolution is computed with Singular \cite{Singular}. Since the resolution is linear the mapping $\alpha$ is given by a $12\times 11$ matrix of linear forms and by the Hilbert-Burch Theorem (see \cite[Theorem 201.5]{Eis95}) the $12$ generators of $\mathcal{J}$ can be chosen as determinants of the maximal order minors of $\alpha$.
\end{proof}

The maximal tangency points were also computed by Edge:
$$\begin{array}{llll}
PK_1=(i:1:-1),\;& PK_2=(-i:1:1),\;& PK_3=(-i:1:-1),\;& PK_4=(i:1:1),\\ 
PK_5=(-1:-i:1),\;& PK_6=(1:-i:1),\;& PK_7=(-1:i:1),\;& PK_8=(1:i:1),\\ 
PK_9=(1:-1:-i),\;& PK_{10}=(1:1:-i),\;& PK_{11}=(1:-1:i),\;& PK_{12}=(1:1:i).
\end{array}$$
Taking the dual lines $PK' = \cald(PK)$ (see Definition \ref{roulleau_D}), we obtain an arrangement of $12$ lines which intersect in $30$ double points and $6$ points of multiplicity $4$, so we obtain the $t$-vector:
$$(30, 0, 6).$$

Moreover, the ideal of the points from $\mathcal{P}_4(PK')$ is a complete intersection ideal generated by $xyz$ and the Fermat conic $x^2+y^2+z^2$.





Interestingly the $30$ double points are contained in the Fermat arrangement of lines determined by linear factors of
$$(x^2-y^2)(y^2-z^2)(z^2-x^2).$$
and there are $5$ such points on each line. Note that points $PK_1, PK_2, \dots, PK_{12}$ also lie on these lines, but they are all different from these $30$ double points.


\begin{proposition}
    Let $\mathcal{K}$ be the saturated ideal of the $30$ points determined by $\mathcal{P}_2(PK')$. Then $\mathcal{K}$ is generated in degrees $6$ and $7$ and its minimal free resolution has the following shape:
    $$0 \longrightarrow S(-9)^3 \longrightarrow S(-6) \oplus S(-7)^3 \longrightarrow S \longrightarrow S/\mathcal{K} \longrightarrow 0.$$
\end{proposition}

\section{Sextactic points on the Fermat quartic}\label{sec4}

A point $P$ on a plane curve $\Gamma$ is called a \emph{sextactic point}, if there exists a curve $C$ of degree $2$, tangent to $\Gamma$ at $P$ with multiplicity at least $6$. We additionally say that $P$ is a \emph{proper sextactic point} if the conic $C$ is reduced.
\begin{remark}
According to our definition inflexion points (flexes) of a curve are sextactic points but not \emph{proper sextactic points}. \\
By the Pl\"ucker Formula there are
$$3d(d-2)$$
inflexion points on a smooth plane curve of degree $d$.
\end{remark}
Sextactic points are of particular interest because in a general point of $\Gamma$ there is a conic intersecting $\Gamma$ with the local intersection multiplicity $5$, see \cite{Cay59} where these points were introduced for the first time. The number of sextactic points on a curve of degree $d$ is equal to
$$3d(4d-9)$$
as computed originally by Cayley \cite{Cay65} and reconfirmed recently as an interesting application of Gromov-Witten invariants by Muratore \cite{Mur24}. In particular for a quartic the number of sextactic points is $84$. Taking away the 24 inflexion points of a quartic we are left with $60$ points. However in the cases studied here, additionally the MTPs need to be taken away. Thus we are left with only $48$ proper sextactic points.
\subsection{Computation for sextactic points}
The exact coordinates of sextactic points on both quartics were already computed in \cite{ALWALEED}, using the method involving symmetric groups. The sizes of orbits of each point under the group action were also computed. Here, we compute these coordinates using a different method, involving the \emph{second hessian}. 

In 1865, Cayley introduced the notion of the second Hessian of a curve $\Gamma$, denoted by $H_2(\Gamma)$. He proved that all sextactic points of the curve $\Gamma$ lie on its second Hessian. He also provided the formula to compute it, which was later corrected by Maugesten and Moe (see \cite{Maugesten_2019} for details).

We used a particularly transparent formula from \cite{SzeSzp24} to compute the second Hessian of the Fermat quartic $F$. We obtain that up to scalar
$$H_2(F) = x^3y^3z^3(x^4 - y^4)(y^4 - z^4)(z^4 - x^4).$$
In particular it splits completely in linear factors. Moreover the coordinate lines contain only inflection points and the MTPs. The ideal of the remaining $48$ sextactic points is thus complete intersection generated by
$$x^4+y^4+z^4\;\mbox{ and }\; (x^4 - y^4)(y^4 - z^4)(z^4 - x^4).$$
This observation opens doors to computing exact coordinates of these points.
\begin{proposition}\label{prop: sextactic_F}
The sextactic points of the Fermat quartic are    
$$\left(1 : \varepsilon^{k-1} : \sqrt[4]{2}\varepsilon^{\ell}\right),\ \left(1 : \sqrt[4]{2}\varepsilon^k : \varepsilon^{\ell - 1}\right),\
\left(1 : {\frac{1}{\sqrt[4]{2}}}\varepsilon^{k} : {\frac{1}{\sqrt[4]{2}}}\varepsilon^{\ell}\right),$$
where $k, \ell \in \{1, 3, 5, 7\}$ and $\varepsilon \in \mathbb{C}$ is the primitive root of unity of order $8$.
\end{proposition}

We also used the formula from \cite{SzeSzp24} to compute the second Hessian of the Komiya-Kuribayashi quartic. We obtain that up to a scalar
$$H_2(K) = xyz(x^2-y^2)(y^2-z^2)(z^2-x^2)Q(x,y,z),$$
where
\begin{align*}
    Q(x,y,z) &= 1056x^{12}-19278x^{10}y^2-75207x^8y^4-111042x^6y^6-75207x^4y^8-19278x^2y^{10}+1056y^{12} \\
    &-19278x^{10}z^2-137198x^8y^2z^2-287194x^6y^4z^2-287194x^4y^6z^2-137198x^2y^8z^2 \\
    &-19278y^{10}z^2-75207x^8z^4-287194x^6y^2z^4-413110x^4y^4z^4-287194x^2y^6z^4-75207y^8z^4\\
    &-111042x^6z^6-287194x^4y^2z^6-287194x^2y^4z^6-111042y^6z^6-75207x^4z^8-137198x^2y^2z^8\\
    &-75207y^4z^8-19278x^2z^{10}-19278y^2z^{10}+1056z^{12}.
\end{align*}

Based on the obtained formula, it is possible to compute the exact coordinates of some of the sextactic points on the Komiya-Kuribayashi quartic. Since every point which lies on both the quartic $K$ and its second Hessian $H_2(K)$ is either a sextactic point or the MTP, we need to find points which lie on both curves and are not on the list $PK_1, \dots, PK_{12}$, introduced earlier.

For instance, one of the factors of $H_2(K)$ are the coordinate lines, while all MTPs have nonzero coordinates.

\begin{proposition}\label{prop: sextactic_K1}
    The points with coordinates
    $$\begin{array}{ll}
        S_1 = \left(0 : 2 : i(\sqrt{5}+1)\right),\ & S_7 = \left(0 : 2 : i(\sqrt{5}-1)\right),\\
        S_2 = \left(0 : 2 : -i(\sqrt{5}+1)\right),\ & S_8 = \left(0 : 2 : -i(\sqrt{5}-1)\right), \\
        S_3 = \left(i(\sqrt{5}+1) : 0 : 2\right),\ & S_9 = \left(i(\sqrt{5}-1) : 0 : 2\right), \\
        S_4 = \left(-i(\sqrt{5}+1) : 0 : 2\right),\ & S_{10} = \left(-i(\sqrt{5}-1) : 0 : 2\right), \\
        S_5 = \left(2 : i(\sqrt{5}+1) : 0\right),\ & S_{11} = \left(2 : i(\sqrt{5}-1) : 0\right), \\
        S_6 = \left(2 : -i(\sqrt{5}+1) : 0\right),\ & S_{12} = \left(2 : -i(\sqrt{5}-1) : 0\right)
    \end{array}$$
    are sextactic points on the Komiya-Kuribayashi quartic.
\end{proposition}

We can find more sextactic points lying on lines obtained from the factors $(x^2 - y^2)$, $(y^2 - z^2)$ and $(z^2 - x^2)$. These lines contain all the MTPs as well as some additional points.

\begin{proposition}\label{prop: sextactic_K2}
    The points with coordinates
    $$\begin{array}{ll}
        S_{13} = (1 : 1 : i\sqrt{5}),\ & S_{19} = (1 : -1 : i\sqrt{5}),\\
        S_{14} = (1 : 1 : -i\sqrt{5}),\ & S_{20} = (1 : -1 : -i\sqrt{5}), \\
        S_{15} = (1 : i\sqrt{5} : 1),\ & S_{21} = (-1 : i\sqrt{5} : 1), \\
        S_{16} = (1 : -i\sqrt{5} : 1),\ & S_{22} = (-1 : -i\sqrt{5} : 1), \\
        S_{17} = (i\sqrt{5} : 1 : 1),\ & S_{23} = (i\sqrt{5} : 1 : -1), \\
        S_{18} = (-i\sqrt{5} : 1 : 1),\ & S_{24} = (-i\sqrt{5} : 1 : -1)
    \end{array}$$
    are sextactic points on the Komiya-Kuribayashi quartic.
\end{proposition}

The remaining $24$ sextactic points lie on the curve $Q(x,y,z) = 0$ mentioned above. It is worth noting that all MTPs on the Komiya-Kuribayashi quartic lie on this curve as well. We can therefore saturate those points from the ideal generated by $Q$ and the Komiya-Kuribayashi quartic. By doing so, we obtained the following result.

\begin{proposition}
    The remaining $24$ sextactic points $S_{25}, \dots, S_{48}$ of the Komiya-Kuribayashi quartic are the complete intersection of that quartic and a curve of degree $6$, given by
    \[90x^2y^4+124x^2y^2z^2+90x^2z^4+45y^6+135y^4z^2+135y^2z^4+45z^6 = 0.\]
    
\end{proposition}

\section{Conics tangent to both quartics at sextactic points}\label{sec5}

In this section, we consider conics which are tangent to the quartic at two sextactic points, with multiplicity at least $4$ in both points. We consider the sextactic points on both the Fermat quartic $F$ and the Komiya-Kuribayashi quartic $K$. We want to give explicit formulas for such conics as well as study their intersection patterns.

In order to find the formulas of such conics, we use the algorithm presented in \cite{MERTA}.

We begin with $48$ sextactic points on $F$, coordinates of which are given in Proposition \ref{prop: sextactic_F}. We obtain the following result.

\begin{theorem}
There are exactly $24$ conics, which are tangent to the Fermat quartic $F$ at two sextactic points, with multiplicity $4$ in both points. These conics are given by the following formulas:
\begin{align*}
    x^2 \pm xy + y^2 \pm \frac{i}{\sqrt{2}}z^2 &= 0, & x^2 \pm ixy - y^2 \pm \frac{i}{\sqrt{2}}z^2 &= 0, \\
    y^2 \pm yz + z^2 \pm \frac{i}{\sqrt{2}}x^2 &= 0, & y^2 \pm iyz - z^2 \pm \frac{i}{\sqrt{2}}x^2 &= 0, \\
    z^2 \pm zx + x^2 \pm \frac{i}{\sqrt{2}}y^2 &= 0, & z^2 \pm izx - x^2 \pm \frac{i}{\sqrt{2}}y^2 &= 0.
\end{align*}
\end{theorem}

It is worth noting that for each sextactic point on $F$ there is only one such a conic passing through it and each conic passes through exactly $2$ sextactic points, therefore connecting these points in $24$ pairs.

Our next step was to figure out how these $24$ conics intersect.

\begin{proposition}
    The conics tangent to $F$ at two sextactic points with multiplicity $4$ in both points intersect in exactly $960$ points: $912$ ordinary double points, $24$ points where two conics are tangent with multiplicity $2$ (tacnodes) and $24$ points where $4$ conics intersect (ordinary quadruple points).
\end{proposition}

We are going to study these two obtained sets of $24$ points in a bit more detail. The points in which two conics are tangent have the following coordinates:
$$\begin{array}{cccc}
    (\pm 2 : i\sqrt{3} + 1 : 0), & (\pm 2 : 0 : i\sqrt{3} + 1), & (0 : \pm 2 : i\sqrt{3} + 1), \\
    (\pm 2i : i\sqrt{3} + 1 : 0), & (\pm 2i : 0 : i\sqrt{3} + 1), & (0 : \pm 2i : i\sqrt{3} + 1), \\
    (i\sqrt{3} + 1 : \pm 2 : 0), & (i\sqrt{3} + 1 : 0 : \pm 2), & (0 : i\sqrt{3} + 1 : \pm 2), \\
    (i\sqrt{3} + 1 : \pm 2i : 0), & (i\sqrt{3} + 1 : 0 : \pm 2i), & (0 : i\sqrt{3} + 1 : \pm 2i), \\  
\end{array}$$

Note that those points lie on the coordinate lines $x = 0$, $y = 0$ and $z = 0$ and there are $8$ points on each line. Moreover, it can be shown that these $24$ points lie on a curve of degree $8$, given by
$$x^8 + y^8 + z^8 + x^4y^4 + y^4z^4 + z^4x^4 = 0,$$
and these points are actually a complete intersection of the above curve and $xyz = 0$.

It turns out that the other set of $24$ points (where $4$ conics intersect) has similar properties. The coordinates of those points are shown below. As before, $\varepsilon$ denotes the primitive $8$th root of unity.
$$\begin{array}{cccc}
    (\pm 1 : \varepsilon\sqrt[4]{2} : 0), & (\pm 1 : 0 : \varepsilon\sqrt[4]{2}), & (0 : \pm 1 : \varepsilon\sqrt[4]{2}), \\
    (\pm i : \varepsilon\sqrt[4]{2} : 0), & (\pm i : 0 : \varepsilon\sqrt[4]{2}), & (0 : \pm i : \varepsilon\sqrt[4]{2}), \\
    (\varepsilon\sqrt[4]{2} : \pm 1 : 0), & (\varepsilon\sqrt[4]{2} : 0 : \pm 1), & (0 : \varepsilon\sqrt[4]{2} : \pm 1), \\
    (\varepsilon\sqrt[4]{2} : \pm i : 0), & (\varepsilon\sqrt[4]{2} : 0 : \pm i), & (0 : \varepsilon\sqrt[4]{2} : \pm i), \\
\end{array}$$

These points lie on the coordinate lines as well, $8$ points on each one. Moreover, they are also a complete intersection, but this time $xyz = 0$ intersects with a curve given by
$$x^8 + y^8 + z^8 + \frac{5}{2}(x^4y^4 + y^4z^4 + z^4x^4) = 0.$$

Now we are going to focus on sextactic points on the Komiya-Kuribayashi quartic. We are going to use only the first $24$ points, coordinates of which we managed to compute. Just like in Proposition \ref{prop: sextactic_K1} and Proposition \ref{prop: sextactic_K2}, we divide them into two groups of $12$ points. In this case, we obtain the following results.

\begin{theorem}
There are exactly $12$ conics, which are tangent to the Komiya-Kuribayashi quartic $K$ at two sextactic points from the set $\{S_1, S_2, \dots, S_{12}\}$, with multiplicity $4$ in both points. These conics are given by the following formulas:
\begin{align*}
    3x^2 + 2y^2 + (-1)^k 2iyz + 2z^2 &= 0, & 6x^2 + \left(5 + (-1)^k \sqrt{5}\right)y^2 + \left(5 + (-1)^{k+1} \sqrt{5}\right)z^2 &= 0, \\
    3y^2 + 2z^2 + (-1)^k 2ixz + 2x^2 &= 0, & 6y^2 + \left(5 + (-1)^k \sqrt{5}\right)z^2 + \left(5 + (-1)^{k+1} \sqrt{5}\right)x^2 &= 0, \\
    3z^2 + 2x^2 + (-1)^k 2ixy + 2y^2 &= 0, & 6z^2 + \left(5 + (-1)^k \sqrt{5}\right)x^2 + \left(5 + (-1)^{k+1} \sqrt{5}\right)y^2 &= 0,
\end{align*}
for $k \in \{0,1\}.$
\end{theorem}

Therefore, there are two conics passing through each point from the set $\{S_1, S_2, \dots, S_{12}\}$ and naturally these conics are tangent in these points with multiplicity $4$. It turns out that apart from the points $\{S_1, S_2, \dots, S_{12}\}$ these conics intersect only in ordinary double points.
 
\begin{theorem}
There are exactly $6$ conics, which are tangent to the Komiya-Kuribayashi quartic $K$ at two sextactic points from the set $\{S_{13}, S_{14}, \dots, S_{24}\}$, with multiplicity $4$ in both points. These conics are given by the following formulas:
\begin{align*}
    4x^2 + 11y^2 \pm 2yz + 11z^2 &= 0, \\
    4y^2 + 11z^2 \pm 2xz + 11x^2 &= 0, \\
    4z^2 + 11x^2 \pm 2xy + 11y^2 &= 0.
\end{align*}
\end{theorem}

This time there is only one conic passing through each point from the set $\{S_{13}, S_{14}, \dots, S_{24}\}$, hence just like in the case of sextactic points on the Fermat quartic $F$, these points are connected into pairs. Moreover, it can be checked that these $6$ conics only intersect in $60$ ordinary double points.

\section*{Appendix}

The code used to obtain the results presented in this paper is available at the following link: 
\\
\url{https://github.com/LukaszM111/Singular-code---quartics}.

All procedures can be run using Singular \cite{Singular}.

\bibliographystyle{abbrv}
\bibliography{master}

\end{document}